\documentstyle[twoside,amsmath,amssymb,12pt]{article}

\renewcommand{\title}[1]{\leftline{\Large\bf #1}\par\medskip}
\renewcommand{\author}[1]{\medskip{\large #1}\par\medskip}
\newcommand{\address}[1]{\noindent{\sl #1}}

\allowdisplaybreaks
\newcommand{\rom}[1]{{\rm #1}}

\begin{document}

\setcounter{page}{1}
\setcounter{section}{0}
\thispagestyle{empty}

\newtheorem{definition}{Definition}
\newtheorem{remark}{Remark}
\newtheorem{proposition}{Proposition}
\newtheorem{theorem}{Theorem}
\newtheorem{corollary}{Corollary}

\title{LAPLACE OPERATORS AND DIFFUSIONS}
\title{IN TANGENT BUNDLES OVER POISSON}
\title{SPACES}

\author{SERGIO ALBEVERIO, ALEXEI DALETSKII,}
\author{AND EUGENE LYTVYNOV}

\begin{abstract}
Spaces of differential forms over configuration spaces with Poisson measures
are constructed. The corresponding Laplacians (of Bochner and de Rham type)
on 1-forms and associated semigroups are considered. Their probabilistic
interpretation is given.
\end{abstract}

\section{Introduction}

Stochastic differential geometry of infinite-dimensional manifolds has been
a very active topic of research in recent times. One of the important and
intriguing problems discussed concerns the construction of spaces of
differential forms over such manifolds and the study of the corresponding
Laplace operators and associated (stochastic) cohomologies. A central role
in this framework is played by the concept of the Dirichlet operator of a
differentiable measure, which is actually an infinite-dimensional
generalization of the Laplace--Beltrami operator on functions, respectively
 the
Laplace--Witten--de Rham operator on differential forms. The study of the
latter operator and the associated semigroup on finite-dimensional manifolds
was the subject of many works, and it leads to deep results on the border of
stochastic analysis, differential geometry and topology, and mathematical
physics, see, e.g., \cite{E3}, \cite{CFKSi}, \cite{EL}. Dirichlet forms on
Clifford algebras were considered in \cite{Gr}. In an infinite-dimensional
situation, such questions were discussed in the flat case in \cite{Ar1},
\cite{Ar2}, \cite{ArM}, \cite{AK}. A regularized heat semigroup on
differential forms over the infinite-dimensional torus was studied in \cite
{LeBe}. A study of such questions on general infinite product manifolds was
given in \cite{ADK1}, \cite{ADK2}. The case of loop spaces was considered in
\cite{JL}, \cite{LRo}.

At the same time, there is a growing interest in geometry and analysis on
Poisson spaces, i.e., on  spaces of locally finite configurations in
non-compact manifolds, equipped with Poisson measures. In \cite{AKR-1}, \cite
{AKR0}, \cite{AKR1}, an approach to these spaces as to infinite-dimensional
manifolds was initiated. This approach was motivated by the connection of
such spaces with the theory of representations of  diffeomorphism groups,
see  \cite{GGPS}, \cite{VGG}, \cite{I} (these references and \cite
{AKR1}, \cite{AKR3} also contain discussion of relations with quantum
physics). We refer the reader to \cite{AKR2}, \cite{AKR3}, \cite{Ro}, and
references therein for further discussion of analysis on Poisson spaces and
applications.

In the present work, we develop this point of view. We define spaces of
differential forms over Poisson spaces. Next, we define and study Laplace
operators acting in the spaces of 1-forms. We show, in particular, that the
corresponding de Rham Laplacian can be expressed in terms of the Dirichlet
operator on functions on the Poisson space and the Witten Laplacian on the
initial manifold associated with the intensity of the corresponding Poisson
measure. We give a probabilistic interpretation and investigate some
properties of the associated semigroups.
Let us remark that the study of Laplacians on $n$-forms by our methods is also
possible, but it leads to more complicated constructions. It will be given in a
forthcoming paper.
The main general aim of our
approach is to develop a framework which extends to Poisson spaces (as
infinite-dimensional manifolds) the finite-dimensional Hodge--de Rham theory.

A different approach to the construction of differential forms and related
objects over Poisson spaces, based on the ``transfer principle'' from
Wiener spaces, is proposed in \cite{Pr2}, see also \cite{PPr} and \cite{Pr}.

\section{Differential forms over configuration spaces}

The aim of this section is to define differential forms over configuration spaces
(as infinite-dimensional manifold).
First, we recall some known facts and definitions concerning ``manifold-like''
 structures
and functional calculus on these spaces.

\subsection{Functional calculus on configuration spaces}

 Our presentation in this subsection is
based on \cite{AKR1}, however for later use in the present paper we give
 a different description of some objects and results occurring in \cite
{AKR1}.

Let $X$ be a complete, connected, oriented, $C^\infty $ (non-compact) Riemannian
manifold of dimension $d$. We denote by $\langle \bullet ,\bullet \rangle _x$
the corresponding inner product in the tangent space $T_xX$ to $X$ at a
point $x\in X$. The associated norm will be denoted by $|\bullet |_x$. Let
also $\nabla ^X$ stand for the gradient on $X$.

The configuration space $\Gamma _X$ over $X$ is defined as the set of all
locally finite subsets (configurations) in $X$:
\begin{equation}\notag
\Gamma _X:=\left\{ \,\gamma \subset X\mid |\gamma \cap \Lambda |<\infty
\text{ for each compact }\Lambda \subset X\,\right\} .
\end{equation}
Here, $|A|$ denotes the cardinality of the set $A$.

We can identify any $\gamma \in \Gamma _X$ with the positive integer-valued
Radon measure
\begin{equation}\notag
\sum_{x\in \gamma }\varepsilon _x\subset {\cal M}(X),
\end{equation}
where $\varepsilon _x$ is the Dirac measure with mass at $x$, $\sum_{x\in
\varnothing }\varepsilon _x:=$zero measure, and ${\cal M}(X)$ denotes the set
of all positive Radon measures on the Borel $\sigma $-algebra ${\cal B}(X)$.
The space $\Gamma _X$ is endowed with the relative topology as a subset of
the space ${\cal M}(X)$ with the vague topology, i.e., the weakest topology
on $\Gamma _X$ such that all maps
\begin{equation}\notag
\Gamma _X\ni \gamma \mapsto \langle f,\gamma \rangle :=\int_Xf(x)\,\gamma
(dx)\equiv \sum_{x\in \gamma }f(x)
\end{equation}
are continuous. Here, $f\in C_0(X)$($:=$the set of all continuous functions
on $X$ with compact support). Let ${\cal B}(\Gamma _X)$ denote the
corresponding Borel $\sigma $-algebra.

Following \cite{AKR1}, we define the tangent space to $\Gamma _X$ at a point
$\gamma $ as the Hilbert space
\begin{equation}\notag
T_\gamma \Gamma _X:=L^2(X\to TX;d\gamma ),
\end{equation}
or equivalently
\begin{equation}
T_\gamma \Gamma _X=\bigoplus_{x\in \gamma }T_xX.  \label{tg-sp1}
\end{equation}
The scalar product and the norm in $T_\gamma \Gamma _X$ will be denoted by $%
\langle \bullet ,\bullet \rangle _\gamma $ and $\left\| \bullet \right\| _\gamma $%
, respectively. Thus, each $V(\gamma )\in T_\gamma \Gamma _X$ has the form $%
V(\gamma )=(V(\gamma )_x)_{x\in \gamma }$, where $V(\gamma )_x\in T_xX$, and
\begin{equation}\notag
\| V(\gamma )\| _\gamma ^2=\sum_{x\in \gamma }|V(\gamma )_x|_x^2.
\end{equation}

Let $\gamma \in \Gamma _X$ and $x\in \gamma $. By ${\cal O}_{\gamma ,x}$ we
will denote an arbitrary open neighborhood of $x$ in $X$ such that the
intersection of the closure of ${\cal O}_{\gamma ,x}$ in $X$ with $\gamma
\setminus \{x\}$ is the empty set. For any fixed finite subconfiguration $%
\left\{ x_1,\dots ,x_k\right\} \subset \gamma $, we will always consider open
neighborhoods ${\cal O}_{\gamma ,x_1},\dots ,{\cal O}_{\gamma ,x_k}$ with
disjoint closures.

Now, for a measurable function $F\colon\Gamma _X\to {\Bbb R}$, $\gamma \in \Gamma
_X$, and $\left\{ x_1,\dots ,x_k\right\} \subset \gamma $, we define a
function $F_{x_1,\dots ,x_k}(\gamma ,\bullet )\colon{\cal O}_{\gamma ,x_1}\times
\dots \times {\cal O}_{\gamma ,x_k}\to {\Bbb R}$ by
\begin{multline*}
{\cal O}_{\gamma ,x_1}\times \dots \times {\cal O}_{\gamma ,x_k} \ni
(y_1,\dots ,y_k)\mapsto F_{x_1,\dots ,x_k}(\gamma ,y_1,\dots ,y_k):= \\
=F((\gamma \setminus \{x_1,\dots ,x_k\})\cup \{y_1,\dots ,y_k\})\in {\Bbb R%
}.
\end{multline*}
Since we will be interested only in the local behavior of the function $%
F_{x_1,\dots ,x_k}(\gamma ,\bullet )$ around the point $(x_1,\dots ,x_k)$, we
will not write explicitly which neighborhoods ${\cal O}_{\gamma ,x_i}$ we
use.

\begin{definition}
\label{def2.0}\rom{We say that a function $F:\Gamma _X\to {\Bbb R}^1$ is
differentiable at $\gamma \in \Gamma _X$ if for each $x\in \gamma $ the
function $F_x(\gamma ,\cdot )$ is differentiable at $x$ and
\[
\nabla ^\Gamma F(\gamma )=(\nabla ^\Gamma F(\gamma )_x)_{x\in \gamma }\in
T_\gamma \Gamma _X,
\]
where
\begin{equation}\notag
\nabla ^\Gamma F(\gamma )_x:=\nabla ^XF_x(\gamma ,x).
\end{equation}
}
\end{definition}

{We will call $\nabla ^\Gamma F(\gamma )$ the gradient of $F$ at $\gamma $. }

For a function $F$ differentiable at $\gamma $ and a vector $V(\gamma )\in
T_\gamma \Gamma _X$, the directional derivative of $F$ at the point $\gamma $
along $V(\gamma )$ is defined by
\begin{equation}\notag
\nabla _V^\Gamma F(\gamma ):=\langle \nabla ^\Gamma F(\gamma ),V(\gamma
)\rangle _\gamma .
\end{equation}

In what follows, we will also use the shorthand notation
\begin{equation}
\nabla _x^XF(\gamma ):=\nabla ^XF_x(\gamma ,x),  \label{flick}
\end{equation}
so that $\nabla ^\Gamma F(\gamma )=(\nabla _x^XF(\gamma ))_{x\in \gamma }$.
It is easy to see that the operation $\nabla ^\Gamma $ satisfies the usual
properties of differentiation, including the Leibniz rule.

We define a class ${\cal FC}_{\mathrm b}^\infty (\Gamma _X)$ of smooth cylinder
functions on $\Gamma _X$ as follows:

\begin{definition}
\label{def2.1} \rom{A measurable bounded function $F:\Gamma _X\to {\Bbb R}
^1$ belongs to ${\cal FC}_{\mathrm b}^\infty (\Gamma _X)$ iff:

(i) there exists a compact $\Lambda \subset X$ such that $F(\gamma
)=F(\gamma _\Lambda )$ for all $\gamma \in \Gamma _X$, where $\gamma
_\Lambda :=\gamma \cap \Lambda $;

(ii) for any $\gamma \in \Gamma _X$ and $\left\{ x_1,\dots ,x_k\right\}
\subset \gamma $, $k\in {\Bbb N}$, the function $F_{x_1,\dots ,x_k}(\gamma
,\bullet )$ is infinitely differentiable with derivatives uniformly bounded in
$\gamma $ and $x_1,\dots ,x_k$ (i.e., the majorizing constant depends only
on the order of differentiation but not on the specific choice of $\gamma
\in \Gamma _X$, $k\in {\Bbb N}$, and $\{x_1,\dots ,x_k\}\subset \gamma $).}
\end{definition}

Let us note that, for $F\in {\cal FC}_{\mathrm b}^\infty (\Gamma _X)$, only a finite
number of coordinates of $\nabla ^\Gamma F(\gamma )$ are not equal to zero,
and so $\nabla ^\Gamma F(\gamma )\in T_\gamma \Gamma _X$. Thus, each $F\in
{\cal FC}_{\mathrm b}^\infty (\Gamma _X)$ is differentiable at any point $\gamma \in
\Gamma _X$ in the sense of Definition~\ref{def2.0}.

\begin{remark}
\label{rem2.1}\rom{In \cite{AKR1}, the authors introduced the class ${\cal FC}
_{\mathrm b}^\infty ({\cal D},\Gamma _X)$ of functions on $\Gamma _X$ of the form
\begin{equation}
F(\gamma )=g_F(\left\langle \varphi _1,\gamma \right\rangle ,\dots
,\left\langle \varphi _N,\gamma \right\rangle ),  \label{2.1}
\end{equation}
where $g_F\in C_{\mathrm b}^\infty ({\Bbb R}^N)$ and $\varphi _1,\dots ,\varphi _N\in
{\cal D}:=C_0^\infty (X)$($:=$ the set of all $C^\infty $-functions on $X$
with compact support). Evidently, we have the inclusion
$${\cal FC}_{\mathrm b}^\infty (
{\cal D},\Gamma _X)\subset {\cal FC}_{\mathrm b}^\infty (\Gamma _X),$$ and moreover,
the gradient of $F$ of the form (\ref{2.1}) in the sense of Definition~\ref
{def2.0},
\begin{equation}\notag
\nabla ^\Gamma F(\gamma )_x=\sum_{i=1}^N\frac{\partial g_F}{\partial s_i}%
(\langle \varphi _1,\gamma \rangle ,\dots ,\langle \varphi
_N,\gamma \rangle )\nabla ^X\varphi _i(x),
\end{equation}
coincides with the gradient of this function in the sense of \cite{AKR1}. }
\end{remark}

\subsection{Tensor bundles and cylinder forms over configuration
\\ spaces}

Our next aim is to introduce differential forms on $\Gamma _X$.

Vector fields and first order differential forms on $\Gamma _X$ will be
identified with sections of the bundle $T\Gamma _X.$ Higher order
differential forms will be identified with  sections of tensor
bundles $\wedge ^n(T\Gamma _X)$ with fibers
\begin{equation}\notag
\wedge ^n(T_\gamma \Gamma _X)%
=%
\wedge ^n(L^2(X\rightarrow TX;\gamma )),
\end{equation}
where $\wedge ^n({\cal H})$ (or ${\cal H}^{\wedge n}$) stands for the $n$-th
antisymmetric tensor power of a Hilbert space ${\cal H}$. In what follows,
we will use different representations of this space. Because of (\ref{tg-sp1}%
), we have
\begin{equation}
\wedge ^n(T_\gamma \Gamma _X)=\wedge ^n\bigg( \bigoplus_{x\in \gamma
}T_xX\bigg) .  \label{tang-n}
\end{equation}

Let us introduce the factor space $X^n/S_n$, where $S_n$ is the permutation
group of $\{1,\dots ,n\}$ which naturally acts on $X^n$:
\begin{equation}\notag
\sigma (x_1,\dots ,x_n)=(x_{\sigma (1)},\dots ,x_{\sigma (n)}),\qquad\sigma \in
S_n.
\end{equation}
The space $X^n/S_n$ consists of equivalence classes $[x_1,\dots ,x_n]$ and
we will denote by $[x_1,\dots ,x_n]_d$ an equivalence class $[x_1,\dots
,x_n] $ such that the equality $x_{i_1}=x_{i_2}=\dots =x_{i_k}$ can hold
only for $k\le d$ points. (In other words, any equivalence class $[x_1,\dots
,x_n]$ is a multiple configuration in $X$, while $[x_1,\dots ,x_n]_d$ is a
multiple configuration with multiplicity of points $\le d$.) We will omit the
lower index $d$ in the case where $n\le d$. In what follows, instead of
writing $[x_1,\dots ,x_n]_d:\{x_1,\dots ,x_n\}\subset \gamma $, we will use
the shortened notation $[x_1,\dots ,x_n]_d\subset \gamma $, though $%
[x_1,\dots ,x_n]_d$ is not, of course, a set. We then have from (\ref{tang-n}%
):
\begin{equation}
\wedge ^n(T_\gamma \Gamma _X)=\bigoplus_{[x_1,\dots ,x_n]_d\,\subset \gamma
}T_{x_1}X\wedge T_{x_2}X\wedge \dots \wedge T_{x_n}X,  \label{tang-n0}
\end{equation}
since for each $\sigma\in S_n$ the spaces $T_{x_1}X\wedge T_{x_2}X\wedge\dots\wedge T_{x_n}X$
and $T_{x_{\sigma(1)}}X\wedge T_{x_{\sigma(2)}}X\wedge\dots\wedge T_{x_{\sigma(n)}}X$ coincide.

Thus, under a differential form $\omega $ of order $n$, $n\in {\Bbb N}$,
over $\Gamma _X,$ we will understand the mapping
\begin{equation}\notag
\Gamma _X\ni \gamma \mapsto \omega (\gamma )\in \wedge ^n(T_\gamma \Gamma
_X).
\end{equation}
We denote by $\omega (\gamma )_{[x_1,\dots ,x_n]_d}$ the corresponding
component in the decomposition (\ref{tang-n0}).

In particular, in the case $n=1$, a 1-form $V$ over $\Gamma _X$ is given by
the mapping
\begin{equation}\notag
\Gamma _X\ni \gamma \mapsto V(\gamma )=(V(\gamma )_x)_{x\in \gamma }\in
T_\gamma \Gamma _X.
\end{equation}

For  fixed $\gamma \in \Gamma _X$ and $x\in \gamma ,$ we consider the
mapping
\begin{equation}\notag
{\cal O}_{\gamma ,x}\ni y\mapsto \omega _x(\gamma ,y)%
:=%
\omega (\gamma _y)\in \wedge ^n(T_{\gamma _y}\Gamma _X),
\end{equation}
where $\gamma _y%
:=%
(\gamma \setminus \{x\})\cup \{y\},$ which is a section of the Hilbert
bundle
\begin{equation}
\wedge ^n(T_{\gamma _y}\Gamma _X)\mapsto y\in {\cal O}_{\gamma ,x}
\label{bund1}
\end{equation}
over ${\cal O}_{\gamma ,x}.$ The Levi--Civita connection on $TX$ generates in
a natural way a connection on this bundle. We denote by $\nabla _{\gamma
,x}^X$ the corresponding covariant derivative, and use the notation
\begin{equation}\notag
\nabla _x^X\omega (\gamma )%
:=%
\nabla _{\gamma ,x}^X\omega _x(\gamma ,x)\in T_xX\otimes \left( \wedge
^n(T_\gamma \Gamma _X)\right)
\end{equation}
if the section $\omega _x(\gamma ,\cdot )$ is differentiable at $x$.
Analogously, we denote by $\Delta _x^X$ the corresponding Bochner Laplacian
associated with the volume measure $m$ on ${\cal O}_{\gamma ,x}$ (see
subsec.~3.2 where the notion of Bochner Laplacian is
recalled).

Similarly, for a fixed $\gamma \in \Gamma _X$ and $\left\{
x_1,\dots,x_n\right\} \subset \gamma $, we define a mapping
\begin{multline*}
{\cal O}_{\gamma ,x_1}\times \dots\times {\cal O}_{\gamma ,x_n} \ni
(y_1,\dots,y_n)\mapsto \omega _{x_1,\dots,x_n}(\gamma ,y_1,\dots,y_n)
:= \\
=
\omega (\gamma _{y_1,\dots,y_n}) \in \wedge ^n(T_{\gamma
_{y_1,\dots,y_n}}\Gamma _X),
\end{multline*}
where $\gamma _{y_1,\dots,y_n}%
:=%
(\gamma \setminus \{x_1,\dots,x_n\})\cup \{y_1,\dots,y_n\}$, which is a section
of the Hilbert bundle
\begin{equation}
\wedge ^n(T_{\gamma _{y_1,\dots,y_n}}\Gamma _X)\mapsto \left(
y_1,\dots,y_n\right) \in {\cal O}_{\gamma ,x_1}\times \dots\times {\cal O}%
_{\gamma ,x_n}  \label{bund-n}
\end{equation}
over ${\cal O}_{\gamma ,x_1}\times\dots\times {\cal O}_{\gamma ,x_n}.$

Let us remark that, for any $\eta \subset \gamma $, the space $\wedge
^n(T_\eta \Gamma _X)$ can be identified in a natural way with a subspace
of $\wedge ^n(T_\gamma \Gamma _X)$. In this sense, we will use
expressions of the type $\omega (\gamma )=\omega (\eta )$ without additional
explanations.

A set ${\cal F}\Omega ^n$ of smooth cylinder $n$-forms over $\Gamma _X$ will
be defined as follows.

\begin{definition}
\label{def2.2}\rom{${\cal F}\Omega ^n$ is the set of $n$-forms $\omega $ over $%
\Gamma _X$ which satisfy the following conditions:

(i) there exists a compact $\Lambda =\Lambda (\omega )\subset X$ such that $%
\omega (\gamma )=\omega (\gamma _\Lambda )$;

(ii) for each $\gamma \in \Gamma _X$ and $\left\{ x_1,...,x_n\right\}
\subset \gamma $, the section $\omega _{x_1,\dots,x_n}(\gamma ,\bullet )$ of the
bundle (\ref{bund-n}) is infinitely differentiable at $(x_1,\dots,x_n),$ and
bounded together with the derivatives uniformly in $\gamma $.
}\end{definition}

\begin{remark}
\label{form-fin}\rom{For each $\omega \in {\cal F}\Omega ^n$, $\gamma \in \Gamma
_X$, and any open bounded $\Lambda \supset \Lambda (\omega )$, we can define
the form $\omega _{\Lambda ,\gamma }$ on ${\cal O}_{\gamma ,x_1}\times \dots
\times {\cal O}_{\gamma ,x_n}$ by
\begin{equation}
\omega _{\Lambda ,\gamma }(y_1,\dots,y_n)=\operatorname{Proj}_{\wedge
^n(T_{y_1}X\oplus \dots \oplus T_{y_n}X)}\omega (\gamma \setminus
\{x_1,\dots,x_n\}\cup \{y_1,\dots,y_n\}),  \label{cyl-form}
\end{equation}
where $\{x_1,\dots,x_n\}=\gamma \cap \Lambda $. The item (ii) of Definition~%
\ref{def2.2} is obviously equivalent to the assumption  $\omega
_{\Lambda ,\gamma }$ to be smooth and bounded together with the derivatives
uniformly in $\gamma $ (for some $\Lambda $ and consequently for any $%
\Lambda \supset \Lambda (\omega )$).}
\end{remark}

\begin{definition}
\label{def2.3}\rom{We define the covariant derivative $\nabla ^\Gamma \omega $ of
the form $\omega \in {\cal F}\Omega ^n$ as the mapping
\begin{equation}\notag
\Gamma _X\ni \gamma \mapsto \nabla ^\Gamma \omega (\gamma )%
:=%
(\nabla _x^X\omega (\gamma ))_{x\in \gamma }\in T_\gamma \Gamma _X\otimes
\left( \wedge ^n(T_\gamma \Gamma _X)\right)
\end{equation}
if for all $\gamma\in\Gamma_X$ and $x\in\gamma$ the form $\omega_x(\gamma,\bullet)$
is differentiable at $x$ and the $\nabla^\Gamma\omega(\gamma)$ just defined
indeed belongs to $T_\gamma \Gamma _X\otimes
\left( \wedge ^n(T_\gamma \Gamma _X)\right) $.
}
\end{definition}

\begin{remark}\rom{
For each $\omega\in{\cal F\Omega}^n$, the covariant derivative $\nabla^\Gamma\omega$
exists, and moreover only a finite number of the coordinates
$\nabla^\Gamma\omega(\gamma)_{x,[x_1,\dots,x_n]_d}$ in the decomposition
$$T_\gamma\Gamma_X\otimes\big(
\wedge^n(T_\gamma\Gamma_X)
\big)=\bigoplus_{x\in\gamma,\,[x_1,\dots,x_n]_d\subset\gamma}
T_xX\otimes(T_{x_1}X\wedge\dots\wedge T_{x_n}X)
$$
are not equal to zero.}
\end{remark}

\begin{proposition}
\label{prop2.1} For arbitrary $\omega ^{(1)},\omega ^{(2)}\in {\cal F}\Omega
^n$, we have
\begin{multline*}\nabla ^\Gamma \langle \omega ^{(1)}(\gamma ),\omega ^{(2)}(\gamma
)\rangle _{\wedge ^n(T_\gamma \Gamma _X)} = \\
=\langle\nabla ^\Gamma \omega ^{(1)}(\gamma ),\omega ^{(2)}(\gamma
)\rangle _{\wedge ^n(T_\gamma \Gamma _X)}+\langle \omega
^{(1)}(\gamma ),\nabla ^\Gamma \omega ^{(2)}(\gamma )\rangle _{\wedge
^n(T_\gamma \Gamma _X)}.
\end{multline*}
\end{proposition}

\noindent
{\it Proof}. We have, for any fixed $\gamma \in \Gamma _X$ and $x\in \gamma $,%
\begin{gather*}
\nabla _x^X\langle \omega ^{(1)}(\gamma ),\omega ^{(2)}(\gamma )\rangle
_{\wedge ^n(T_\gamma \Gamma _X)}=\nabla _x^X\langle \omega _x^{(1)}(\gamma
,x),\omega _x^{(2)}(\gamma ,x)\rangle _{\wedge ^n(T_\gamma \Gamma _X)}
 \\
=\langle \nabla _x^X\omega ^{(1)}(\gamma ),\omega ^{(2)}(\gamma )\rangle
_{\wedge ^n(T_\gamma \Gamma _X)}+\langle \omega ^{(1)}(\gamma ),\nabla
_x^X\omega ^{(2)}(\gamma )\rangle _{\wedge ^n(T_\gamma \Gamma _X)},
\end{gather*}
because of the usual properties of the covariant derivative $\nabla _x^X$.
\quad$\blacksquare $

\subsection{Square integrable forms}

In this subsection, we will consider spaces of forms over the configuration space
$\Gamma_X$ which are square integrable with respect to a Poisson measure.

Let $m$ be the volume measure on $X$, let $\rho \colon X\to {\Bbb R}$ be a
measurable function such that $\rho >0$ $m$-a.e., and $\rho ^{1/2}\in
H_{\mathrm loc}^{1,2}(X)$, and define the measure $\sigma (dx):=\rho (x)\,m(dx)$.
Here, $H_{\mathrm loc}^{1,2}(X)$ denotes the local Sobolev space of order 1 in $%
L_{\mathrm loc}^2(X;m)$. Then, $\sigma $ is a non-atomic Radon measure on $X$.

Let $\pi _\sigma $ stand for the Poisson measure on $\Gamma _X$ with
intensity $\sigma $. This measure is characterized by its Fourier transform
\begin{equation}\notag
\int_{\Gamma _X}e^{i\langle f,\gamma \rangle }\,\pi _\sigma (d\gamma )=\exp
\int_X(e^{if(x)}-1)\,\sigma (dx),\qquad f\in C_0(X).
\end{equation}
Let $F\in L^1(\Gamma _X;\pi _\sigma )$ be cylindrical, that is, there exits
a compact $\Lambda \subset X$ such that $F(\gamma )=F(\gamma _\Lambda )$.
Then, one has the following formula, which we will use many times:
\begin{equation}
\int_{\Gamma _X}F(\gamma )\,\pi _\sigma (d\gamma )=e^{-\sigma (\Lambda
)}\sum_{n=0}^\infty \frac 1{n!}\int_{\Lambda ^n}F(\{x_1,\dots
,x_n\})\,\sigma (dx_1)\dotsm\sigma (dx_n).  \label{3.1}
\end{equation}
Since the measure $\sigma $ is non-atomic, the sets $\{(x_1,\dots ,x_n)\in
\Lambda ^n:x_i=x_j\}$, $i,j=1,\dots ,n$, $i\ne j$, have zero $\sigma
(dx_1)\dotsm\sigma (dx_n)$ measure, and therefore the expression on the right
hand side of (\ref{3.1}) is well-defined.

We define on the set ${\cal F}\Omega ^n$ the $L^2$-scalar product with
respect to the Poisson measure:
\begin{equation}
(\omega ^{(1)},\omega ^{(2)})_{L_{\pi _\sigma }^2\Omega ^n}%
:=%
\int_{\Gamma _X}\langle \omega ^{(1)}(\gamma ),\omega ^{(2)}(\gamma )\rangle
_{\wedge ^nT_\gamma \Gamma _X}\,\pi _\sigma (d\gamma ).  \label{4.1}
\end{equation}

As easily seen, for each $\omega \in {\cal F}\Omega ^n$, the function $%
\langle \omega (\gamma ),\omega (\gamma )\rangle _{\wedge ^nT_\gamma \Gamma
_X}$ is polynomially bounded on $\Gamma _X$, and therefore it belongs to all
$L^p(\Gamma _X;\pi _\sigma )$, $p\ge 1$. Moreover, $(\omega ,\omega
)_{L_{\pi _\sigma }^2\Omega ^n}>0$ if $\omega $ is not identically zero.
Hence, we can define the Hilbert space
\begin{equation}\notag
L_{\pi _\sigma }^2\Omega ^n%
:=%
L^2(\Gamma _X\to \wedge ^nT\Gamma _X;\pi _\sigma )
\end{equation}
as the closure of ${\cal F}\Omega ^n$ in the norm generated by the scalar
product (\ref{4.1}).

From now on, we consider the case of 1-forms only and suppose that $\dim
X\ge 2$. We give another description of the spaces $L_{\pi _\sigma }^2\Omega
^1.$ Let us recall the following well-known result (Mecke identity, see e.g.
\cite{KMM}):
\begin{equation}
\int_{\Gamma _X}\int_Xf(\gamma ,x)\,\gamma (dx)\,\pi (d\gamma )=\int_{\Gamma
_X}\int_Xf(\gamma \cup \{x\},x)\,\sigma (dx)\,\pi _\sigma (d\gamma )
\label{mecke}
\end{equation}
for any measurable bounded $f:\Gamma _X\times X\rightarrow {\Bbb R}^1$.

Let us introduce the notations $$L_\sigma ^2\Omega ^1(X)%
:=%
L^2(X\rightarrow TX;\sigma ),\qquad
L_{\pi _\sigma }^2(\Gamma _X)%
:=
L^2(\Gamma _X\rightarrow {\Bbb R}^1;\pi _\sigma ).$$

\begin{proposition}\label{kkkk}
The space $L_{\pi _\sigma }^2\Omega ^1$ is isomorphic to the space $L_{\pi_\sigma}
^2(\Gamma _X)\otimes L_\sigma ^2\Omega ^1(X)$ with the isomorphism $I^1$
given by the formula
\begin{equation}
( I^1V) (\gamma ,x)%
:=%
V(\gamma \cup \{x\})_x,\qquad\gamma \in \Gamma _X,\;x\in X.  \label{isom}
\end{equation}
\end{proposition}

\noindent
{\it Proof}. Let us specify the scalar product of two cylinder 1-forms $%
V,W\in {\cal F}\Omega ^1.$ We have:
\begin{align*}
( W,V) _{L_{\pi _\sigma }^2\Omega ^1} &=\int_{\Gamma
_X}\left\langle W(\gamma ),V(\gamma )\right\rangle _\gamma\, \pi _\sigma
(d\gamma )  \\
&=\int_{\Gamma _X}\int_X\left\langle W(\gamma )_x,V(\gamma )_x\right\rangle
_\gamma \,\gamma (dx)\,\pi _\sigma (d\gamma ) \\
&=\int_{\Gamma _X}\int_X\left\langle W(\gamma \cup \{x\})_x,V(\gamma \cup
\{x\})_x\right\rangle _\gamma\, \gamma (dx)\,\pi _\sigma (d\gamma ),
\end{align*}
because $\gamma \cup \{x\}=\gamma $ for $x\in \gamma .$ The application of
the Mecke identity to the function
\begin{equation}\notag
f(\gamma ,x)=\left\langle V(\gamma \cup \{x\})_x,W(\gamma \cup
\{x\})_x\right\rangle _\gamma
\end{equation}
shows that
\begin{equation}\notag
\left( W,V\right) _{L_{\pi _\sigma }^2\Omega ^1}=\int_{\Gamma
_X}\int_X\left\langle V(\gamma \cup \{x\})_x,W(\gamma \cup
\{x\})_x\right\rangle _\gamma \,\sigma (dx)\,\pi _\sigma (d\gamma ).
\end{equation}
The space ${\cal F}\Omega ^1$ is, by definition, dense in $L_{\pi _\sigma
}^2\Omega ^1,$ and so it remains only to show that $I^1({\cal F}\Omega ^1)$
is dense in $L_{\pi _\sigma }^2(\Gamma _X)\otimes L_\sigma ^2\Omega ^1(X).$

For $F\in {\cal FC}_{\mathrm b}^\infty (\Gamma _X),$ and $\nu \in \Omega _0^1(X)$ (the
set of smooth 1-forms on $X$ with compact support), we define a form $%
V(\gamma )=(V(\gamma )_x)_{x\in \gamma }$ by setting
\begin{equation}
V(\gamma )_x%
:=%
F(\gamma \setminus \{x\})\nu (x).  \label{isom1}
\end{equation}
Evidently, we have $V\in {\cal F}\Omega ^1,$ and
\begin{equation}
\left( I^1V\right) (\gamma ,x)=F(\gamma )\nu (x)  \label{isom2}
\end{equation}
for each $\gamma $ and any $x\notin \gamma $. Since each $\gamma \in \Gamma
_X$ is a subset of $X$ of zero $m$ measure, we conclude from (\ref{isom2})
that
\begin{equation}
I^1V=F\otimes \nu .  \label{isom3}
\end{equation}
Noting that the linear span of such $F\otimes \nu $ is dense in $L_{\pi
_\sigma }^2(\Gamma _X)\otimes L_\sigma ^2\Omega ^1(X)$, we obtain the
result. \quad $\blacksquare $

In what follows, we will denote by ${\cal D}\Omega ^1$ the linear span of
forms $V$ defined by (\ref{isom1}). As we
 already noticed in the proof of Proposition~\ref{kkkk}, ${\cal D}\Omega
^1\subset {\cal F}\Omega ^1$ and is dense in $L_{\pi _\sigma }^2\Omega ^1$.

\begin{corollary}
\label{fock}We have the unitary isomorphism
\begin{equation}\notag
{\cal I}:L_{\pi _\sigma }^2\Omega ^1\rightarrow \operatorname{Exp}\left( L^2(X;\sigma
)\right) \otimes L_\sigma ^2\Omega ^1(X)
\end{equation}
given by
\begin{equation}\notag
{\cal I}=\left( U\otimes {\bf 1}\right) I^1,
\end{equation}
where $U$ is the unitary isomorphism between the Poisson space $L_{\pi
_\sigma }^2(\Gamma _X)$ and the symmetric Fock space $\operatorname{Exp}\left( L^2(X;\sigma
)\right) $\rom, see e\rom.g\rom{.\ \cite{AKR1}.}
\end{corollary}

\section{Dirichlet operators on differential forms over configuration
spaces\label{dodf}}

In this section, we introduce Dirichlet operators associated with the
Poisson measure on $\Gamma _X$ which act in the space $L^2_{\pi_\sigma}
\Omega^1$. These operators generalize the notions of Bochner and de
Rham--Witten Laplacians on finite dimensional manifolds.

In the  two first subsections, we recall some known facts and definitions
concerning Dirichlet operators of Poisson measures on configuration spaces
and Laplace operators on differential forms over finite-dimensional manifolds.

\subsection{The intrinsic Dirichlet operator on functions}

In this subsection, we recall some theorems from \cite{AKR1}
which concern the  intrinsic Dirichlet operator in the space $L^2_{\pi_\sigma}
(\Gamma_X)$, to be used later.

Let us recall that the logarithmic derivative of the measure $\sigma $ is
given by the vector field
\begin{equation}\notag
X\ni x\mapsto \beta _\sigma (x):=\frac{\nabla ^X\rho (x)}{\rho (x)}\in T_xX
\end{equation}
(where as usual $\beta _\sigma :=0$ on $\{\rho =0\}$). We wish now to define
the notion of logarithmic derivative of the Poisson measure, and for this we
need a generalization of the notion of vector field.

For each $\gamma \in \Gamma _X$, consider the triple
\begin{equation}\notag
T_{\gamma ,\,\infty }\Gamma _X\supset T_\gamma \Gamma _X\supset T_{\gamma
,0}\Gamma _X.
\end{equation}
Here, $T_{\gamma ,0}\Gamma _X$ consists of all finite sequences from $%
T_\gamma \Gamma _X$, and $T_{\gamma ,\,\infty }\Gamma _X%
:=%
\left( T_{\gamma ,0}\Gamma _X\right) ^{\prime }$ is the dual space, which
consists of all sequences $V(\gamma )=(V(\gamma )_x)_{x\in \gamma }$, where $%
V(\gamma )_x\in T_xX$. The pairing between any $V(\gamma )\in T_{\gamma
,\,\infty }\Gamma _X$ and $v(\gamma )\in T_{\gamma ,0}\Gamma _X$ with
respect to the zero space $T_\gamma \Gamma _x$ is given by
\begin{equation}\notag
\langle V(\gamma ),v(\gamma )\rangle _\gamma =\sum_{x\in \gamma }\langle
V(\gamma )_x,v(\gamma )_x\rangle _x
\end{equation}
(the series is, in fact, finite). From now on, under a vector field over $%
\Gamma _X$ we will understand mappings of the form $\Gamma _X\ni \gamma
\mapsto V(\gamma )\in T_{\gamma ,\infty }\Gamma _X$.

The logarithmic derivative of the Poisson measure $\pi _\sigma $ is defined
as the vector field
\begin{equation}
\Gamma _X\ni \gamma \mapsto B_{\pi _\sigma }(\gamma )=(\beta _\sigma
(x))_{x\in \gamma }\in T_{\gamma ,\infty }\Gamma _X
\end{equation}
(i.e., the logarithmic derivative of the Poisson measure is the lifting of
the logarithmic derivative of the underlying measure).

The following theorem is a version of Theorem~3.1 in \cite{AKR1} (for more
general classes of functions and vector fields).

\begin{theorem}[Integration by parts formula on the Poisson space]
\label{th-ibp}$\text{}$\newline
For arbitrary $F^{(1)},F^{(2)}\in {\cal FC}_{\mathrm b}^\infty
(\Gamma _X)$ and a smooth cylinder vector field $V\in
{\cal FV}(\Gamma_X)$ $(:=
{\cal F}\Omega^1)$\rom, we have
\begin{gather*}
\int_{\Gamma _X}\nabla _V^\Gamma F^{(1)}(\gamma )F^{(2)}(\gamma )\,\pi
_\sigma (d\gamma )=-\int_{\Gamma _X}F^{(1)}(\gamma )\nabla _V^\Gamma
F^{(2)}(\gamma )\,\pi _\sigma (d\gamma )  \\
\text{}
-\int_{\Gamma _X}F^{(1)}(\gamma )F^{(2)}(\gamma )\big[ \left\langle B_{\pi
_\sigma }(\gamma ),V(\gamma )\right\rangle _\gamma +\operatorname{div}^\Gamma
V(\gamma )\big] \,\pi _\sigma (d\gamma ),
\end{gather*}
where the divergence $\operatorname{div}^\Gamma V(\gamma )$ of the vector field $V$
is given by
\begin{gather*}
\operatorname{div}V(\gamma )=\sum_{x\in \gamma }\operatorname{div}_x^XV(\gamma
)=\langle \operatorname{div}_{\bullet }^XV(\gamma ),\gamma \rangle ,
\\
\operatorname{div}_x^XV(\gamma )%
:=%
\operatorname{div}^XV_x(\gamma ,x),\qquad x\in \gamma ,
\end{gather*}
$\operatorname{div}^X$ denoting the divergence on $X$ with respect to the volume
measure $m.$
\end{theorem}

\noindent
{\it Proof}. The theorem follows from formula (\ref{3.1}) and the usual
integration by parts formula on the space $L^2(\Lambda ^n,\sigma ^{\otimes
n})$ (see also the proof of Theorem~\ref{th4.1} below).\quad $\blacksquare $
\vspace{2mm}

Following \cite{AKR1}, we consider the intrinsic pre-Dirichlet form on the
Poisson space
\begin{equation}
{\cal E}_{\pi _\sigma }(F^{(1)},F^{(2)})=\int_{\Gamma _X}\langle \nabla
^\Gamma F^{(1)}(\gamma ),\nabla ^\Gamma F^{(2)}(\gamma )\rangle _\gamma
\,\pi _\sigma (d\gamma )  \label{3.2}
\end{equation}
with domain $D({\cal E}_{\pi _\sigma }):={\cal FC}_{\mathrm b}^\infty (\Gamma _X)$. By
using the fact that the measure $\pi _\sigma $ has all moments finite, one
can show that the expression (\ref{3.2}) is well-defined.

Let $H_\sigma $ denote the Dirichlet operator in the space $L^2(X;\sigma )$
associated to the pre-Dirichlet form
\begin{equation}\notag
{\cal E}_\sigma (\varphi ,\psi )=\int_X\langle \nabla ^X\varphi (x),\nabla
^X\psi (x)\rangle _x\,\sigma (dx),\qquad \varphi ,\psi \in {\cal D}.
\end{equation}
This operator acts as follows:
\begin{equation}\notag
H_\sigma \varphi (x)=-\Delta ^X\varphi (x)-\langle \beta _\sigma (x),\nabla
^X\varphi (x)\rangle _x,\qquad \varphi\in{\cal D},
\end{equation}
where $\Delta ^X:=\operatorname{div}^X\nabla ^X$ is the Laplace--Beltrami operator on
$X$.

Then, by using Theorem~\ref{th-ibp}, one gets
\begin{equation}
{\cal E}_{\pi _\sigma }(F^{(1)},F^{(2)})=\int_{\Gamma _X}H_{\pi _\sigma
}F^{(1)}(\gamma )F^{(2)}(\gamma )\,\pi _\sigma (d\gamma ),\qquad
F^{(1)},F^{(2)}\in {\cal FC}_{\mathrm b}^\infty (\Gamma _X).  \label{sukk}
\end{equation}
Here, the intrinsic Dirichlet operator $H_{\pi _\sigma }$ is given by
\begin{align}
H_{\pi _\sigma }F(\gamma )
:&=%
\sum_{x\in \gamma }H_{\sigma ,x}F(\gamma )\equiv \langle H_{\sigma ,\bullet
}F(\gamma ),\gamma \rangle , \notag  \\
H_{\sigma ,x}F(\gamma ) %
:&=%
H_\sigma F_x(\gamma ,x),\qquad x\in \gamma , \label{dir-op1}
\end{align}
so that the operator $H_{\pi _\sigma }$ is the lifting to $L^1(\Gamma _X;\pi
_\sigma )$ of the operator $H_\sigma $ in $L^2(X;\sigma )$.

Upon (\ref{sukk}), the pre-Dirichlet form ${\cal E}_{\pi _\sigma }$ is
closable, and we preserve the notation for the closure of this form.

\begin{theorem}
\label{th3.2}\rom{\cite{AKR1}} Suppose that $(H_\sigma ,{\cal D})$ is essentially
self-adjoint on $L^2(X;\sigma )$\rom. Then\rom,
the operator $H_{\pi _\sigma }$ is
essentially self-adjoint on ${\cal FC}_{\mathrm b}^\infty (\Gamma _X).$
\end{theorem}

\begin{remark}
\label{proof}\rom{This theorem was proved in \cite{AKR1}, Theorem~5.3. (We have
already mentioned in Remark~\ref{rem2.1} that the inclusion ${\cal FC}%
_{\mathrm b}^\infty ({\cal D},\Gamma _X)\subset {\cal FC}_{\mathrm b}^\infty (\Gamma _X)$
holds.) We would like to stress that this result is based on the theorem
which says that the image of the operator $H_{\pi _\sigma }$ under the
isomorphism $U$ between the Poisson space and the  Fock space $\operatorname{Exp%
}\left( L^2(X;\sigma )\right) $over $L^2(X;\sigma )$ is the differential
second quantization $d\operatorname{Exp}H_\sigma $ of the operator $H_\sigma $.
}\end{remark}

\begin{remark}\rom{
In what follows, we will always assume that the conditions of the theorem
are satisfied. It is true e.g.\ in the case where $\| \beta _\sigma
\| _{TX}\in L_{\mathrm loc}^p(X;\sigma )$ for some $p>\dim \,X,$ see \cite
{AKR1}.
}\end{remark}

Finally, we mention the important fact that the diffusion process which is
properly associated with the Dirichlet form $({\cal E}_{\pi _\sigma },D(%
{\cal E}_{\pi _\sigma }))$ is the usual independent infinite particle
process (or distorted Brownian motion), cf. \cite{AKR1}.

\subsection{Laplacians on differential forms over
finite-dimensional\\
manifolds}

We recall now some facts on the Bochner and de Rham--Witten Laplacians
on differential forms
over
a finite-dimensional manifold.

Let $M$ be a Riemannian manifold equipped with the measure $\mu (dx)=e^{\phi
(x)}dx,$ $dx$ being the volume measure and $\phi $  a $C_{\mathrm b}^2$-function
on $M$. We consider a Hilbert bundle
\begin{equation}\notag
{\cal H}_x\mapsto x\in M
\end{equation}
over $M$ equipped with a smooth connection, and denote by $\nabla $ the
corresponding covariant derivative in the spaces of sections of this bundle.
Let $L^2(M\rightarrow {\cal H};\mu )$ be the space of $\mu $-square
integrable sections. The operator
\begin{equation}\notag
H_\mu ^B%
:=%
\nabla _\mu ^{*}\nabla
\end{equation}
in $L^2(M\rightarrow {\cal H};\mu ),$ where $\nabla _\mu ^{*}$ is the
adjoint of $\nabla $, will be called the Bochner Laplacian associated with
the measure $\mu $. Differentiability of $\mu $ implies that $\nabla _\mu
^{*}\nabla $ is a uniquely defined self-adjoint operator. One can easily
write the corresponding differential expression on the space of twice
differentiable sections. In the case where $\phi \equiv 0$ and ${\cal H}%
_x=\wedge ^n(T_xM)$, we obtain the classical Bochner Laplacian on differential
forms (see \cite{CFKSi}).

Now, let $d$ be the exterior differential in spaces of differential forms
over $M.$ The operator
\begin{equation}\notag
H_\mu ^R%
:=%
d_\mu ^{*}d+dd_\mu ^{*}
\end{equation}
acting in the space of $\mu $-square integrable forms, where $d_\mu ^{*}$ is
the adjoint of $d$, will be called the de Rham Laplacian associated with the
measure $\mu $ (or the Witten Laplacian associated with $\phi $, see \cite
{CFKSi}).

The relation of the Bochner and de Rham--Witten Laplacians on differential
forms is given by the Weitzenb\"{o}ck formula, which in the case of 1-forms
has the following form (see\cite{CFKSi}, \cite{E3}):
\begin{equation}\notag
H_\mu ^Ru(x)=H_\mu ^Bu(x)+R_\mu (x)u(x),
\end{equation}
where
\begin{equation}
R_\mu (x)%
:=%
R(x)-\nabla ^X\beta _\mu (x).  \label{weitz}
\end{equation}
Here, $R(x)\in {\cal L}(T_xM)$ is the usual Weitzenb\"{o}ck correction term:
\begin{equation}\notag
R(x)%
:=%
\sum_{i,j=1}^{\dim M}\operatorname{Ric}_{ij}(x)a_i^{*}a_j,
\end{equation}
where $\operatorname{Ric}$ is the Ricci tensor on $X$, and $a_i^{*}$ and $a_j$ are
the creation and annihilation operators, respectively.

\subsection{Bochner Laplacian on 1-forms over the Poisson space}

Let us consider the pre-Dirichlet form
\begin{equation}
{\cal E}_{\pi _\sigma }^B(W^{(1)},W^{(2)})=\int_{\Gamma _X}\langle \nabla
^\Gamma W^{(1)}(\gamma ),\nabla ^\Gamma W^{(2)}(\gamma )\rangle _{T_\gamma
\Gamma _X\otimes T_\gamma \Gamma _X}\,\pi _\sigma (d\gamma ),  \label{4.2}
\end{equation}
where $W^{(1)},W^{(2)}\in {\cal F}\Omega ^1$. Again using the fact that $\pi
_\sigma $ has finite moments, one shows that the function under the sign of
integral in (\ref{4.2}) is integrable with respect to $\pi _\sigma $.

\begin{theorem}
\label{th4.1} For any $W^{(1)},W^{(2)}\in {\cal F}\Omega ^1$\rom, we have
\begin{equation}\notag
{\cal E}_{\pi _\sigma }^B(W^{(1)},W^{(2)})=\int_{\Gamma _X}\left\langle
H_{\pi _\sigma }^BW^{(1)}(\gamma ),W^{(2)}(\gamma )\right\rangle _{T_\gamma
\Gamma _X}\,\pi _\sigma (d\gamma ),
\end{equation}
where $H_{\pi _\sigma }^B$ is the operator in the space $L_{\pi _\sigma
}^2\Omega ^1$ given by
\begin{equation}
H_{\pi _\sigma }^BW=-\Delta ^\Gamma W-\left\langle \nabla ^\Gamma W,B_{\pi
_\sigma }(\gamma )\right\rangle _\gamma ,\qquad W\in {\cal F}\Omega ^1.
\label{boch1}
\end{equation}
Here\rom,
\begin{equation}
\Delta ^\Gamma W(\gamma )%
:=%
\sum_{x\in \gamma }\Delta _x^XW(\gamma )\equiv \left\langle \Delta _{\bullet
}^\Gamma W(\gamma ),\gamma \right\rangle ,  \label{boch2}
\end{equation}
where $\Delta _x^X$ is the Bochner Laplacian of the bundle $T_{\gamma
_y}\Gamma _X\mapsto y\in {\cal O}_{\gamma ,x}$ with the volume measure\rom.
\end{theorem}

\noindent
{\it Proof}. Let us fix $W^{(1)},W^{(2)}\in {\cal F}\Omega ^1$. Let $\Lambda
$ be an open bounded set in $X$ such that $\Lambda (W^{(1)})\subset \Lambda$,
$\Lambda (W^{(2)})\subset \Lambda $. Then, by using (\ref{3.1}),
\begin{gather*}
\int_{\Gamma _X}\langle \nabla ^\Gamma W^{(1)}(\gamma ),\nabla ^\Gamma
W^{(2)}\rangle _{T_\gamma \Gamma _X\otimes T_\gamma \Gamma _X}\,\pi _\sigma
(d\gamma ) = \\
=e^{-\sigma (\Lambda )}\sum_{k=0}^\infty \frac 1{k!}\int_{\Lambda
^k}\sum_{i=1}^k\langle \nabla _{x_i}^XW^{(1)}(\{x_1,\dots
,x_k\}) ,  \\
 \nabla _{x_i}^XW^{(2)}(\{x_1,\dots ,x_k\})\rangle
_{T_{x_i}X\otimes T_{\{x_1,\dots ,x_k\}}\Gamma _X}\,\sigma (dx_1)\dotsm\sigma
(dx_k)   \\
=e^{-\sigma (\Lambda )}\sum_{k=0}^\infty \frac
1{k!}\sum_{i=1}^k\int_{\Lambda ^k}\langle \Delta
_{x_i}^XW^{(1)}(\{x_1,\dots ,x_k\}) \\\text{}+
\langle \nabla _{x_i}^XW^{(1)}(\{x_1,\dots ,x_k\}),\beta _\sigma
(x_i)\rangle _{T_{x_i}X},   \\
 W^{(2)}(\{x_1,\dots ,x_k\})\rangle _{T_{\{x_1,\dots
,x_k\}}\Gamma _X}\,\sigma (dx_1)\dotsm\sigma (dx_k)  \\
=\int_{\Gamma _X}\langle H_{\pi _\sigma }^BW^{(1)}(\gamma ),W^{(2)}(\gamma
)\rangle _{T_\gamma \Gamma _X}\,\pi _\sigma (d\gamma ). \quad \blacksquare
\end{gather*}

\begin{remark}
\label{rem4.1}\rom{We can rewrite the action of the operator $H_{\pi _\sigma }^B$
in the two following forms:

\begin{enumerate}
\item[1)]  We have from (\ref{boch1}) and (\ref{boch2}) that
\begin{equation}
H_{\pi _\sigma }^BW(\gamma )=\sum_{x\in \gamma }H_{\sigma ,x}^BW(\gamma
)\equiv \left\langle H_{\sigma ,\bullet }^BW(\gamma ),\gamma \right\rangle
,\qquad W(\gamma )\in {\cal F}\Omega ^1,  \label{blo1}
\end{equation}
where
\begin{equation}
H_{\sigma ,x}^BW(\gamma ):=-\Delta _x^XW(\gamma )-\left\langle \nabla
_x^XW(\gamma ),\beta _\sigma (x)\right\rangle _x.  \label{blo2}
\end{equation}
Thus, the operator $H_{\sigma ,x}^B$ is the lifting of the Bochner Laplacian
on $X$ with the measure $\sigma .$

\item[2)]  As  easily seen, the operator $H_{\pi _\sigma }^B$ preserves
the space ${\cal F}\Omega ^1$, and we can always take $\Lambda (H_{\pi
_\sigma }^BW)=\Lambda (W)$. Then for any open bounded $\Lambda \supset
\Lambda (W)$
\begin{equation}
(H_{\pi _\sigma }^BW)_{\Lambda ,\gamma }=H_{\sigma ,\Lambda \cap \gamma
}^BW_{\Lambda ,\gamma },  \label{cyl-boch}
\end{equation}
where $H_{\sigma ,\Lambda \cap \gamma }^B$ is the Bochner Laplacian of the
manifold $X^{\Lambda \cap \gamma }%
:=%
\times _{x\in \Lambda \cap \gamma }X_x$, $X_x\equiv X$, with the product
measure $\sigma ^{\Lambda \cap \gamma }%
:=%
\otimes _{x\in \Lambda \cap \gamma }\sigma _x$, $\sigma _x\equiv \sigma $
(cf. (\ref{cyl-form})).
\end{enumerate}
}\end{remark}

It follows from Theorem~\ref{th4.1} that the pre-Dirichlet form ${\cal E}%
_{\pi _\sigma }^B$ is closable in the space $L_{\pi _\sigma }^2\Omega ^1$.
The generator of its closure (being actually the
Friedrichs extension of the operator $H_{\pi _\sigma }^B$,
for which we will use the same notation) will be called the Bochner
Laplacian on 1-forms over $\Gamma _X$ corresponding to the Poisson measure $%
\pi _\sigma $.

For operators $A$ and $B$ acting in Hilbert spaces ${\cal H}$ and ${\cal K}$,
respectively, we introduce the operator $A\boxplus B$ in ${\cal H\otimes K}$
by
\begin{equation}\notag
A\boxplus B%
:=%
A\otimes {\bf 1}+{\bf 1}\otimes B.
\end{equation}

\begin{proposition}
\rom{1)}  On ${\cal D}\Omega ^1$ we have
\begin{equation}
I^1\,H_{\pi _\sigma }^B=\left( H_{\pi _\sigma }\boxplus H_\sigma ^B\right)
\,I^1.  \label{dec-gen0}
\end{equation}

\rom{2)}  ${\cal D}\Omega ^1$ is a domain of essential self-adjointness of $%
H_{\pi _\sigma }^B.$

\end{proposition}

\noindent
{\it Proof}.
1) Let $W\in {\cal D}\Omega ^1.$ Then, for some $F\in {\cal FC}_{\mathrm b}^\infty
(\Gamma _X)$, $\omega\in \Omega _0^1(X)$, and any $\gamma \in \Gamma _X$, $x,z\in
\gamma $, $y\in {\cal O}_x$, we have
\begin{equation}\notag
W_x(\gamma ,y)_z=\begin{cases}
F((\gamma \setminus\{x,z\})\cup \{y\})\omega(z),& z\ne y, \\
 F(\gamma \setminus \{x\})\omega(y),&z=y\end{cases}
\end{equation}
Thus
\begin{equation}\notag
H_{\sigma ,x}^BW(\gamma )_z=\begin{cases}
H_{\sigma ,x}F(\gamma \backslash \{z\})\omega(z),&z\ne x, \\
 F(\gamma \setminus \{z\})H_\sigma ^B\omega(z),& z=x.
\end{cases}
\end{equation}
Formula (\ref{dec-gen0}) follows now from (\ref{blo1}) and (\ref{isom2}).

2) The statement follows from (\ref{dec-gen0}) and the essential
self-adjointness of $H_{\pi _\sigma }$ on ${\cal FC}_{\mathrm b}^\infty (\Gamma _X)$
(Theorem~\ref{th3.2}) and $H_\sigma ^B$ on $\Omega _0^1(X)$ (the latter fact
can be shown by standard methods similar to \cite{E2}, \cite{E3}) by the
theory of operators admitting separation of variables \cite[Ch.6]{B}.\quad $
\blacksquare $\vspace{2mm}

We give also a Fock space representation of the operator $H_{\pi _\sigma
}^B $. Corollary~\ref{fock} implies the following

\begin{corollary}
We have
\begin{equation}\notag
{\cal I}H_{\pi _\sigma }^B{\cal I}^{-1}=d\operatorname{Exp}H_\sigma \boxplus H_\sigma ^B,
\end{equation}
cf\rom.\ Remark~\rom{ \ref{proof}.}
\end{corollary}

\subsection{De Rham Laplacian on 1-forms over the Poisson space}

We define the linear operator
\begin{equation}\notag
d^\Gamma :{\cal F}\Omega ^1\to {\cal F}\Omega ^2
\end{equation}
by
\begin{equation}
(d^\Gamma W)(\gamma ):=\sqrt{2}\,\operatorname{AS}(\nabla _x^XW(\gamma )),  \label{5.1}
\end{equation}
where $\operatorname{AS}:(T_\gamma \Gamma _X)^{\otimes 2}\to (T_\gamma \Gamma _X)^{\wedge
2} $ is the antisymmetrization operator. It follows from this definition
that
\begin{equation}
(d^\Gamma W)(\gamma )=\sum_{x\in \gamma }(d_x^XW)(\gamma ),  \label{5.2}
\end{equation}
where
\begin{align}
(d_x^XW)(\gamma ) %
:&=%
\sum_{y\in \gamma }d^X(W_x(\gamma ,x)_y) \notag \\
&=\sum_{y\in \gamma }\sqrt{2}\,\operatorname{AS}(\nabla ^XW_x(\gamma ,x)_y)  \label{5.3}
\end{align}
with $\operatorname{AS}:T_xX\otimes T_yX\to T_xX\wedge T_yX$ being again the
antisymmetrization. This implies that we have indeed the inclusion $d^\Gamma
W\in {\cal F}\Omega ^2$ for each $W\in {\cal F}\Omega ^1$.

Suppose that, for $W\in{\cal F}\Omega^1$, $\gamma\in\gamma_X$, and $x,y\in\gamma$,
the 1-form $W_x(\gamma,\bullet)_y$ has, in local coordinates on the manifold
$X$, the following form:
\begin{equation}\label{kuku}
W_x(\gamma,\bullet)_y=u(\bullet)h,\qquad  u\colon{\cal O}_{\gamma,x}
\to{\Bbb R},\ h\in T_y.\end{equation}
Then, we have
\begin{equation}\label{krik}
\operatorname{AS}(\nabla^XW_x(\gamma,x)_y)=\nabla^X u(x)\wedge h,
\end{equation}
which, upon \eqref{5.3}, describes the action of $d^X_x$.

Let us consider $d^\Gamma $ as an operator acting from the space $L_{\pi
_\sigma }^2\Omega ^1$ into $L_{\pi _\sigma }^2\Omega ^2$. Analogously to the
proof of Theorem~\ref{th4.1}, we get the following formula for the adjoint
operator $( d_{\pi _\sigma }^\Gamma ) ^{*}$ restricted to ${\cal F%
}\Omega ^2$:
\begin{equation}
( d_{\pi _\sigma }^\Gamma ) ^{*}W(\gamma )=\sum_{x\in \gamma
}( d_{\sigma ,x}^X) ^{*}W(\gamma ),\qquad W\in{\cal F}\Omega^2,  \label{5.4}
\end{equation}
where
\begin{equation}
( d_{\sigma ,x}^X) ^{*}W(\gamma )=\sum_{y\in \gamma }(
d_{\sigma ,x}^X) ^{*}W_x(\gamma ,x)_{[x,y]}.  \label{5.5}
\end{equation}
Suppose that, in local coordinates on the manifold $X$, the form $W_x(\bullet
,\bullet )_{[x,y]}$ has the representation
\begin{equation}
W_x(\gamma ,\bullet)_{[x,y]}=w(\bullet)h_1\wedge h_2,\qquad w\colon
{\cal O}_{\gamma,x}\to{\Bbb R},\ h_1\in T_xX,\ h_2\in T_yX.
\label{loc-form}
\end{equation}
Then, taking to notice \eqref{krik}, one concludes that
\begin{multline}
( d_{\sigma ,x}^X) ^{*}(W_x(\gamma ,x)_{[x,y]})  =
-\frac 1{\sqrt{2}}\left[ \left( \langle \nabla ^Xw(x),h_1\rangle
_x+w(x)\langle \beta _\sigma (x),h_1\rangle _x\right) h_2\right.
\\
\text{}-\left. \delta _{x,y}\left( \langle \nabla ^Xw(x),h_2\rangle
_x+w(x)\langle \beta _\sigma (x),h_2\rangle _x\right) h_1\right]. \label{5.6}
\end{multline}
Here,
\begin{equation}\notag
\delta _{x,y}=
\begin{cases}
1,&\text{if } x=y, \\
0,&\text{otherwise.}
\end{cases}
\end{equation}

In what follows, we will suppose for simplicity that the function $\rho $ is
infinite differentiable on $X$ and $\rho (x)>0$ for all $x\in X$. Then, by (%
\ref{5.4})--(\ref{5.6})
\begin{equation}\notag
( d_{\pi _\sigma }^\Gamma ) ^{*}:{\cal F}\Omega ^2\to {\cal F}%
\Omega ^1.
\end{equation}

We set also
\begin{equation}\label{chc}d^\Gamma:{\cal FC}_{\mathrm b}^\infty(\Gamma_X)\to{\cal F}\Omega^1,\qquad d^\Gamma:=\nabla^\Gamma.
\end{equation}
Evidently, the restriction to ${\cal F}\Omega^1$ of the adjoint of $d^\Gamma$ considered as an operator
acting from $L^2_{\pi_\sigma}(\Gamma_X)$ into $L^2_{\pi_\sigma}\Omega^1$ is given by
\begin{equation}\label{lao}(d_{\pi_\sigma}^\Gamma)^{*}:
{\cal F}\Omega^1\to
{\cal FC}^\infty_{\mathrm b}(\Gamma_X),\qquad (d_{\pi_\sigma}^\Gamma)^{*}V(\gamma)=-\operatorname{div}^\Gamma V(\gamma)-
\langle V(\gamma),B_{\pi_\sigma}(\gamma)\rangle_\gamma.\end{equation}

For $n\in {\Bbb N}$, we define the pre-Dirichlet form ${\cal E}_{\pi _\sigma
}^R$ by
\begin{multline*}
{\cal E}_{\pi _\sigma }^R(W^{(1)},W^{(2)}) %
:=%
\int_{\Gamma _X}\big[ \langle d^\Gamma W^{(1)}(\gamma ),d^\Gamma
W^{(2)}(\gamma )\rangle _{\wedge ^2(T_\gamma \Gamma _X)}
\\
\text{}+\langle ( d_{\pi _\sigma }^\Gamma )
^{*}W^{(1)}(\gamma ),( d_{\pi _\sigma }^\Gamma )
^{*}W^{(2)}(\gamma )\rangle _{T_\gamma \Gamma _X}\big] \,\pi _\sigma
(d\gamma ),
\end{multline*}
where $W^{(1)},W^{(2)}\in D({\cal E}^R_{\pi_\sigma}):={\cal F}\Omega ^1$.

The next theorem follows easily from (\ref{5.1})--(\ref{lao}).

\begin{theorem}
\label{th5.1} For any $W^{(1)},W^{(2)}\in {\cal F}\Omega ^1$\rom, we have
\begin{equation}\notag
{\cal E}_{\pi _\sigma }^R(W^{(1)},W^{(2)})=\int_{\Gamma _X}\langle
H_{\pi _\sigma }^RW(\gamma )^{(1)}(\gamma ),W(\gamma )^{(2)}\rangle
_\gamma \,\pi _\sigma (d\gamma ).
\end{equation}
Here,
\begin{equation}\notag
H_{\pi _\sigma }^R%
:=%
d^\Gamma ( d_{\pi _\sigma }^\Gamma ) ^{*}+( d_{\pi _\sigma
}^\Gamma ) ^{*}d^\Gamma ,\qquad D(H_{\pi _\sigma }^R):={\cal F}\Omega ^1,
\end{equation}
is an operator in the space $L_{\pi _\sigma }^2\Omega ^1$\rom. It can be
represented as follows\rom:
\begin{equation}
H_{\pi _\sigma }^RW(\gamma )=\sum_{x\in \gamma }H_{\sigma ,x}^RW(\gamma
)\equiv \left\langle H_{\sigma ,\bullet }^R\,W(\gamma ),\gamma \right\rangle
,  \label{5.7}
\end{equation}
where
\begin{equation}
H_{\sigma ,x}^R=d_x^X( d_{\sigma ,x}^X) ^{*}+( d_{\sigma
,x}^X) ^{*}d_x^X.  \label{5.8}
\end{equation}
\end{theorem}

From Theorem~\ref{th5.1} we conclude that the pre-Dirichlet form ${\cal E}%
_{\pi _\sigma }^R$ is closable in the space $L_{\pi _\sigma }^2\Omega ^1$.

The generator of its closure (being actually the
Friedrichs extension of the operator $H_{\pi _\sigma }^R$,
for which we will use the same notation)
 will be called the de Rham Laplacian on $\Gamma
_X$ corresponding to the Poisson measure $\pi _\sigma $. By (\ref{5.7}) and (%
\ref{5.8}), $H_{\pi _\sigma }^R$ is the lifting of the de Rham Laplacian on $%
X$ with measure $\sigma $.

\begin{remark}\rom{
Similarly to (\ref{cyl-boch}), the operator $H_{\pi _\sigma }^B$ preserves the
space ${\cal F}\Omega ^1$, and we can always take $\Lambda (H_{\pi _\sigma
}^BW)=\Lambda (W)$. Then for any open bounded $\Lambda \supset \Lambda (W)$,
we have
\begin{equation}
(H_{\pi _\sigma }^RW)_{\Lambda ,\gamma }=H_{\sigma ,\Lambda \cap \gamma
}^RW_{\Lambda ,\gamma },  \label{cyl-der}
\end{equation}
where $H_{\sigma ,\Lambda \cap \gamma }^R$ is the de Rham Laplacian of the
manifold $X^{\Lambda \cap \gamma }$ with the product measure $\sigma
^{\Lambda \cap \gamma }$.}
\end{remark}

\begin{proposition}
\rom{1)}  On ${\cal D}\Omega ^1$ we have
\begin{equation}
I^1\,H_{\pi _\sigma }^R=\left( H_{\pi _\sigma }\boxplus H_\sigma ^R\right)
\,I^1.  \label{dec-gen2}
\end{equation}

\rom{2)}  ${\cal D}\Omega ^1$ is a domain of essential self-adjointness of
$H_{\pi _\sigma }^R.$

\end{proposition}

\noindent {\it Proof}.
1) The proof is similar to that of (\ref{dec-gen0}). It is only
necessary to note that, for a ``constant'' 1-form $W$ such that $W(\gamma
)_x=\nu (x)$, we have evidently $\left( H_{\sigma ,x}^Rw(\gamma )\right)
_x=H_\sigma ^R\omega (x)$.

2) The proof is similar to that of the corresponding statement for the
Bochner Laplacian $H_{\pi _\sigma }^B$.\quad $\blacksquare $

\begin{remark}\rom{
By similar methods, one can define  Bochner and de Rham Laplacians on $n$%
-forms over $\Gamma _X$. An extension to this case of formulas (\ref
{dec-gen0}) and (\ref{dec-gen2}) will have, however, a more complicated form.
}\end{remark}

\subsection{Weitzenb\"{o}ck formula on the Poisson space}

In this section, we will derive a generalization of the Weitzenb\"{o}ck
formula to the case of the Poisson measure on the configuration space. In
other words, we will derive a formula which gives a relation between the
Bochner and de Rham Laplacians. We assume that the Weitzenb\"{o}ck
correction term $R_\sigma (x)\in  {\cal L}(T_xX)$ (cf.\ (\ref{weitz})) is
bounded uniformly in $x\in X$.

Given an operator field
\begin{equation}
X\ni x\mapsto J(x)\in {\cal L}(T_xX)  \label{op-pot}
\end{equation}
on $X$ (with $J(x)$ bounded uniformly in $x\in X$), we define the
``diagonal'' operator field
\begin{equation}
\Gamma _X\ni \gamma \mapsto {\bf J}(\gamma )\in {\cal L}(T_\gamma \Gamma _X),
\label{op-field}
\end{equation}
using the decomposition (\ref{tg-sp1}). Thus, we can define the operator
field ${\bf R}_\sigma (\gamma )$.

\begin{theorem}[Weitzenb\"ock formula on the Poisson space]
We have\rom, for each $W\in {\cal F}\Omega ^1,$%
\begin{equation}
H_{\pi _\sigma }^RW(\gamma )=H_{\pi _\sigma }^BW(\gamma )+{\bf R}_\sigma
(\gamma )W(\gamma ).  \label{6.2}
\end{equation}
\end{theorem}

\noindent
{\it Proof}. Let us fix $W\in {\cal F}\Omega ^1$ and $\gamma \in \Gamma _X$.
Let $\Lambda \subset X$ be an open bounded set such that $\Lambda \supset
\Lambda (W)$ (cf.\ Definition~\ref{def2.2}), and let ${\cal O}_{\gamma
,x_1}\times \dots \times {\cal O}_{\gamma ,x_k}$ and $W_{\Lambda ,\gamma }$
be as in Remark \ref{form-fin}. We have then, according to (\ref{cyl-boch})
and (\ref{cyl-der}),
\begin{align*}
(H_{\pi _\sigma }^BW)_{\Lambda ,\gamma } &=H_{\sigma ,\Lambda \cap \gamma
}^BW_{\Lambda ,\gamma },  \\
(H_{\pi _\sigma }^RW)_{\Lambda ,\gamma } &=H_{\sigma ,\Lambda \cap \gamma
}^RW_{\Lambda ,\gamma },
\end{align*}
and the Weitzenb\"{o}ck formula for the manifold $X^{\Lambda \cap \gamma }$
and the measure $\sigma ^{\Lambda \cap \gamma }$ implies that
\[
H_{\sigma ,\Lambda \cap \gamma }^RW_{\Lambda ,\gamma
}(y_1,\dots,y_k)=H_{\sigma ,\Lambda \cap \gamma }^BW_{\Lambda ,\gamma
}(y_1,\dots,y_k)+R(y_1,\dots,y_k)W_{\Lambda ,\gamma }(y_1,\dots,y_k),
\]
where the correction term $R(y_1,\dots,y_k)\in {\cal L}(T_{(y_1,\dots,y_k)}X^{%
\Lambda \cap \gamma })$ is equal to the restriction of ${\bf R}_\sigma
(\gamma )$ to the space $T_{(y_1,\dots,y_k)}X^{\Lambda \cap \gamma }$
(considered as a subspace of $T_\gamma \Gamma _X$), which
is well-defined because of the ``diagonal'' character of ${\bf R}_\sigma
(\gamma )$. It is now enough to remark that the forms $H_{\pi _\sigma }^RW$
and $H_{\pi _\sigma }^BW$ are completely defined by the corresponding
forms $(H_{\pi _\sigma }^RW)_{\Lambda ,\gamma }$ and  $(H_{\pi _\sigma
}^BW)_{\Lambda ,\gamma }$, respectively\quad $\blacksquare $

We can give also an intrinsic description of the correction term ${\bf R}%
_\sigma ^n(\gamma )$. To this end, for each fixed $\gamma \in \Gamma _X$, we define the
operator $R(\gamma )\colon T_{\gamma,0}\Gamma_X\to T_{\gamma,0}\Gamma_X$ as follows:
\begin{gather}
R(\gamma ) %
:=%
\sum_{x\in \gamma }R(\gamma ,x),%
%
\nonumber \\
 \label{corr}
R(\gamma,x)(V(\gamma)_y):=\delta_{x,y}\sum_{i,j=1}^d\operatorname{Ric}_{ij}(x)e_i\,\langle V(\gamma)_x,e_j\rangle_x,\qquad V(\gamma)\in
T_{\gamma,0}
\Gamma_X.
\end{gather}
Here, $\{e_j\}_{j=1}^d$ is again a fixed orthonormal basis in the space $%
T_xX $ considered as a subspace of $T_\gamma \Gamma _X$.

Next, we note that
\begin{align*}
\nabla ^\Gamma B_{\pi _\sigma }(\gamma ) &=(\nabla _x^XB_{\pi _\sigma
}(\gamma ))_{x\in \gamma }=(\nabla _x^X(B_{\pi _\sigma }(\gamma
)_y))_{x,y\in \gamma }  \nonumber \\
&=(\delta _{x,y}\nabla ^X\beta _\sigma (y))_{x,y\in \gamma }\in (T_{\gamma
,\infty }\Gamma _X)^{\otimes 2}.
\end{align*}
Hence, for any $V(\gamma )\in T_{\gamma ,0}\Gamma _X$,
\begin{align}
\nabla _V^\Gamma B_{\pi _\sigma }(\gamma ) %
:&=%
\langle \nabla ^\Gamma B_{\pi _\sigma }(\gamma ),V(\gamma )\rangle _\gamma
\nonumber \\
&=\bigg( \sum_{y\in \gamma }\delta _{x,y}\langle \nabla ^X\beta _\sigma
(y),V(\gamma )_y\rangle _y\bigg) _{x\in \gamma }  \notag\\
&=\left( \langle \nabla ^X\beta _\sigma (x),V(\gamma )_x\rangle _x\right)
_{x\in \gamma }\in T_{\gamma ,0}\Gamma _X.  \label{sdlog}
\end{align}
Thus, $\nabla ^\Gamma B_{\pi _\sigma }(\gamma )$ determines the linear
operator in $T_{\gamma ,0}\Gamma _X$ given by
\begin{equation}\notag
T_{\gamma ,0}\Gamma _X\ni V(\gamma )\mapsto \nabla ^\Gamma B_{\pi _\sigma
}(\gamma )V(\gamma )%
:=%
\nabla _V^\Gamma B_{\pi _\sigma }(\gamma )\in T_{\gamma ,0}\Gamma _X.
\end{equation}

\begin{proposition} We have
\begin{equation}\notag
{\bf R}_\sigma (\gamma )W(\gamma )=R(\gamma )W(\gamma )-\nabla ^\Gamma
B_{\pi _\sigma }(\gamma )W(\gamma ).
\end{equation}
\end{proposition}

\noindent
{\it Proof}. The proposition is derived
from the definition of $%
{\bf R}_\sigma $ and formulas (\ref{sdlog}) and (\ref{corr}).\quad $\blacksquare $

\section{Probabilistic representations of the Bochner and de Rham Laplacians}

Let $\xi _x(t)$ be the Brownian motion with the drift $\beta _\sigma $ on $%
X $ started at a point $x\in X$. We suppose the following:

\begin{itemize}
\item  for each $x$, the process $\xi _x(t)$ has an infinite life-time;

\item  the semigroup
\begin{equation}\notag
T_0(t)f(x)%
:=%
{\sf E}\,f(\xi _x(t))
\end{equation}
acting in the space of bounded measurable functions on $X$ can be extended
to a strongly continuous semigroup of contractions in $L^2(X;\sigma )$, and
its generator $H_0$ is essentially self-adjoint on the space $\cal D$
(in this case $H_0=-H_\sigma $).
\end{itemize}

It follows from the general theory of stochastic differential equations
that these assumptions are satisfied if e.g.\ $\beta _\sigma \in
C_b^4(X\rightarrow TX).$

We denote by $\xi _\gamma (t)$ the corresponding independent particle
process on $\Gamma _X$ which starts at a point $\gamma ,$%
\begin{equation}\notag
\xi _\gamma (t)=(\xi _x(t))_{x\in \gamma }.
\end{equation}
Let
\begin{equation}
{\bf T}_0(t)F(\gamma )%
:=%
{\sf E}\,F(\xi _\gamma (t))
\end{equation}
be the corresponding semigroup in the space of measurable bounded functions
on $\Gamma _X$. It is shown in \cite{AKR1}, that it can be extended to a
strongly continuous semigroup in $L_{\pi _\sigma }^2(\Gamma _X)$ with the
generator ${\bf H}_0=-H_{\pi _\sigma }$ on ${\cal FC}_{\mathrm b}
^\infty (\Gamma _X)$.

Given the operator field (\ref{op-pot}), which is supposed to be continuous
and symmetric (i.e., $J(x)=J(x)^{*}$), we define the operator
\begin{equation}
{\bf P}_{\xi _\gamma }^J(t):T_{\xi _\gamma (t)}\Gamma _X\rightarrow T_\gamma
\Gamma _X  \label{partr}
\end{equation}
by setting
\begin{equation}\notag
( {\bf P}_{\xi _\gamma }^J(t)V) _x=( P_{\xi _x}^J(t))
^{*}V_{\xi _x(t)},\qquad V\in T_{\xi _\gamma (t)}\Gamma _X,
\end{equation}
where the operator
\[
( P_{\xi _x}^J(t)) ^{*}:T_{\xi _x(t)}X\rightarrow T_xX
\]
is adjoint (w.r.t.\ the Riemannian structure of $X$) of the parallel
translation
\begin{equation}\notag
P_{\xi _x}^J(t):T_xX\rightarrow T_{\xi _x(t)}X
\end{equation}
along $\xi _x(t)$ with potential $J.$ That is, $\eta (t)=P_{\xi _x}^J(t)h$
satisfies the SDE
\begin{equation}
\frac D{dt}\eta (t)=J(\eta (t)),\qquad \eta (0)=h,
\end{equation}
where $\frac D{dt}$ is the covariant differentiation along the paths of the
process $\xi $ (see \cite{E3}). It is known that
\begin{equation}
\| P_{\xi _x}^J(t)\| \le e^{tC},
\notag\end{equation}
where $C$ is the supremum of the spectrum of $J(x)$. This implies obviously
the similar estimate for ${\bf P}_{\xi _\gamma }^J$:
\begin{equation}
\| {\bf P}_{\xi _\gamma }^J(t)\| \le e^{tC}.  \label{est1}
\end{equation}

Let us define a semigroup ${\bf T}_1^{{\bf J}}(t)$ associated with the
process $\xi _\gamma $ and potential ${\bf J}$.

\begin{definition}\rom{
For $V\in {\cal F}\Omega ^1$, we set
\begin{equation}
{\bf T}_1^{{\bf J}}(t)V(\gamma )%
:=%
{\sf E}\,{\bf P}_{\xi _\gamma }^J(t)V(\xi _\gamma (t)).\notag
\end{equation}
}\end{definition}

Let $T_1^J(t)$ be the semigroup acting in $L_\sigma ^2\Omega ^1(X)$ as
\begin{equation}\notag
T_1^J(t)\nu (x)%
:=%
{\sf E\,}P_{\xi _x}^J(t)^{*}\nu(\xi _x(t)).
\end{equation}
The following result describes the structure and properties of the semigroup
${\bf T}_1^{{\bf J}}(t)$.

\begin{proposition}
\label{pnsem}\rom{1)}  ${\bf T}_1^{{\bf J}}(t)$ satisfies the estimate
\begin{equation}
\| {\bf T}_1^{{\bf J}}(t)V(\gamma )\| _\gamma \le e^{tC}{\bf T}%
_0(t)\| V(\gamma )\| _\gamma .  \label{markov}
\end{equation}

\rom{2)}  Under the action of the isomorphism $I^1$\rom,
${\bf T}_1^{{\bf J}}(t)
$ obtains the following form\rom:
\begin{equation}
I^1{\bf T}_1^{{\bf J}}(t)={\bf T}_0(t){\bf \otimes }T_1^J(t)\;I^1.
\label{dec-sem}
\end{equation}

\rom{3)}  ${\bf T}_1^{{\bf J}}(t)$ extends to a strongly continuous
semigroup in $L_{\pi _\sigma }^2\Omega ^1.$

\end{proposition}

\noindent{\it Proof}.
1) The result follows from formula (\ref{est1}).

2) Let $V\in{\cal D}\Omega^1$ be given by \eqref{isom1}.
By the definition of ${\bf T}_1^{{\bf J}}(t)$ and
the construction of the process $\xi _\gamma $, we have
\begin{equation}\notag
{\bf T}_1^{{\bf J}}(t)V(\gamma )={\sf E\,}{\bf P}_{\xi _\gamma }^{{\bf J}%
}(t)V(\xi _\gamma (t))
\end{equation}
and
\begin{align*}
( {\bf T}_1^{{\bf J}}(t)V(\gamma )) _x &={\sf E\,}F(\xi _\gamma
(t)\setminus \left\{ \xi _x(t)\right\} )P_{\xi _x}^{J_1}(t)^{*}\nu(\xi _x(t))
 \\
&={\sf E\,}F(\xi _\gamma (t)\setminus \left\{ \xi _x(t)\right\} )\,{\sf E}%
_{\xi _x} \,P_{\xi _x}^{J_1}(t)^{*}\nu(\xi _x(t)) \\
&={\bf T}_0(t)F(\gamma \setminus \left\{ x\right\} )T_1^J(t)\nu(x),
\end{align*}
${\sf E}_{\xi _x}$ meaning the expectation w.r.t.\ the process $\xi _x(t)$,
from where the result follows.

3) The result follows from the corresponding results for semigroups ${\bf T}%
_0(t){\bf \ }$and $T_1^J(t)$, which are well-known (see \cite{AKR1} resp.
\cite{E3}).\quad $\blacksquare $\vspace{2mm}

 Let ${\bf H}_1^{{\bf J}}$ and  $H_1^J$ be the generators of ${\bf T}_1^{{\bf J}}(t)$
 and $%
T_1^J(t)$, respectively.

Now we give probabilistic representations of the semigroups $T_{\pi _\sigma
}^B(t)$ and $T_{\pi _\sigma }^R(t)$ associated with operators $H_{\pi
_\sigma }^B$ and $H_{\pi _\sigma }^R$, respectively. We set
$J_0=0,$ $J_1(x)=R_\sigma
(x)$ (cf.\ (\ref{weitz})). Let us remark that $P_{\xi _x}^{J_0}(t)\equiv
P_{\xi _x}(t)$ is the parallel translation of 1-forms along the path $\xi
_x,$ and we have $H_1^{J_0}=-H_\sigma ^B$ and $H_1^{J_1}=-H_\sigma ^R$ on $%
\Omega _0^1(X).$ We have the following

\begin{theorem}
\rom{1)}  For $W\in {\cal D}\Omega ^1$\rom, we have
\begin{equation}
H_{\pi _\sigma }^BW=-{\bf H}_1^{{\bf J}_0}W,\qquad H_{\pi _\sigma }^RW=-{\bf H}%
_1^{{\bf J}_1}W.  \label{prrepgen}
\end{equation}

\rom{2)}  As $L^2$-semigroups\rom,
\begin{equation}
T_{\pi _\sigma }^B(t)={\bf T}_1^{{\bf J}_0}(t),\qquad T_{\pi _\sigma }^R(t)={\bf T%
}_1^{{\bf J}_1}(t).  \label{prrepsem}
\end{equation}

\rom{3)}  The semigroups $T_{\pi _\sigma }^B(t)$ and $T_{\pi _\sigma
}^R(t)$ satisfy the estimates\rom:
$$
\| T_{\pi _\sigma }^B(t)V(\gamma )\|_\gamma \le {\bf T}_0(t)\|
V(\gamma )\|_\gamma$$
and
$$
\| T_{\pi _\sigma }^R(t)V(\gamma )\|_\gamma \le e^{tC}{\bf T}%
_0(t)\| V(\gamma )\|_\gamma .
$$

\end{theorem}

\noindent
{\it Proof}. 1)
It follows directly from the decomposition (\ref{dec-sem}) that, on ${\cal D}%
\Omega ^1$, we have
\begin{equation}
I^1\,{\bf H}_1^{{\bf J}}=\left( {\bf H}_0\boxplus H_1^J\right) \,I^1.
\label{dec-gen}
\end{equation}
 Setting $J=J_0$ and $J=J_1$
and comparing (\ref{dec-gen0}) and (%
\ref{dec-gen}), we obtain the result.

2) The statement follows from (\ref{prrepgen}) and the essential
self-adjointness of $H_{\pi_\sigma}^B$ and $H_{\pi_\sigma}^R$ on $%
{\cal D}\Omega ^1$ by applying  Proposition \ref{pnsem}, 3), with $J=J_0$
and $J=J_1$, respectively.

3) The result follows from (\ref{prrepsem}) and (\ref{markov}).\quad
 $\blacksquare $

\section{Acknowledgments}
The first author is very grateful to the organizers
for giving him the possibility to present his results at a most stimulating
conference. It is a great pleasure to thank our friends and colleagues Yuri
Kondratiev, Tobias Kuna, and Michael R\"{o}ckner for their interest in this
work and the joy of collaboration. We are also grateful to V. Liebscher for
a useful discussion. The financial support of SFB 256 and DFG Research
Project AL 214/9-3 is gratefully acknowledged.

\address{Institut f\"{u}r Angewandte Mathematik,
Universit\"{a}t Bonn, Wegelerstr. 6,
D 53115 Bonn; and\\
SFB 256, Univ.~Bonn; and\\
 CERFIM (Locarno); Acc. Arch. (USI); and\\
BiBoS, Univ.\ Bielefeld}\vspace{2mm}

\address{Institut f\"{u}r Angewandte Mathematik,
Universit\"{a}t Bonn, Wegelerstr. 6,
D 53115 Bonn; and\\
SFB 256, Univ.~Bonn; and\\
Institute of Mathematics, Kiev; and\\
BiBoS, Univ.\ Bielefeld}\vspace{2mm}

\address{Institut f\"{u}r Angewandte Mathematik,
Universit\"{a}t Bonn, Wegelerstr. 6,
D 53115 Bonn; and\\
BiBoS, Univ.\ Bielefeld}

\end{document}